\documentclass[12pt,a4paper]{amsart}
\usepackage{amssymb}
\usepackage[mathcal]{eucal}

\theoremstyle{plain}
\newtheorem{teo}{Theorem}[section]
\newtheorem{lem}[teo]{Lemma}
\newtheorem{pro}[teo]{Proposition}
\newtheorem{cor}[teo]{Corollary}

\newtheorem{teoABC}{Theorem}

\newtheorem{proABC}[teoABC]{Proposition}

\theoremstyle{remark}
\newtheorem{eje}[teo]{Example}
\newtheorem{que}{Question}[section]

\numberwithin{equation}{section}

\newcommand{\N}{\mathbb{N}}
\newcommand{\Z}{\mathbb{Z}}
\newcommand{\Q}{\mathbb{Q}}
\renewcommand{\L}{\mathbf{Lie}}
\newcommand{\G}{\mathbf{Grp}}

\renewcommand{\phi}{\varphi}

\DeclareMathOperator{\Lie}{lie}
\DeclareMathOperator{\Grp}{grp}
\DeclareMathOperator{\iso}{iso}

\begin{document}
\title[Analytic pro-$p$ groups of small dimensions]{Analytic pro-$p$
  groups of small dimensions}

\thanks{The first author is supported by the Spanish Ministry of
  Science and Education, grants MTM2004-04665 and MTM2007-65385,
  partly with FEDER funds.}

\author{Jon Gonz\'alez-S\'anchez} \address{Departmento de
  Matem\'aticas, Estad\'{\i}stica y Computaci\'on, Facultad de
  Ciencias, Universidad de Cantabria, Avda. de los Castros, E-39071
  Santander, Spain }

\email{jon.gonzalez@unican.es}

\author{Benjamin Klopsch} \address{Department of Mathematics, Royal
  Holloway, University of London, Egham TW20 0EX, United Kingdom}

\email{Benjamin.Klopsch@rhul.ac.uk}

\begin{abstract}
  According to Lazard, every $p$-adic Lie group contains an open
  pro-$p$ subgroup which is saturable.  This can be regarded as the
  starting point of $p$-adic Lie theory, as one can naturally
  associate to every saturable pro-$p$ group $G$ a Lie lattice $L(G)$
  over the $p$-adic integers.

  Essential features of saturable pro-$p$ groups include that they are
  torsion-free and $p$-adic analytic. In the present paper we prove a
  converse result in small dimensions: every torsion-free $p$-adic
  analytic pro-$p$ group of dimension less than $p$ is saturable.

  This leads to useful consequences and interesting questions. For
  instance, we give an effective classification of $3$-dimensional
  soluble torsion-free $p$-adic analytic pro-$p$ groups for $p >
  3$. Our approach via Lie theory is comparable with the use of
  Lazard's correspondence in the classification of finite $p$-groups
  of small order.
\end{abstract}


\maketitle

\section{Introduction}

In his seminal paper \emph{Groupes analytiques
  $p$-adiques}~\cite{La65} from 1965, Lazard proved that a topological
group is $p$-adic analytic if and only if it contains an open pro-$p$
subgroup which is saturable. Indeed, the class of saturable pro-$p$
groups features prominently in $p$-adic Lie theory.  This is due to
the fact that one can naturally associate to every saturable pro-$p$
group a Lie lattice over the $p$-adic integers $\Z_p$. In fact, Lazard
established an isomorphism between the category of saturable pro-$p$
groups and the category of saturable $\Z_p$-Lie lattices.

In the 1980s, Lubotzky and Mann reinterpreted the group-theoretic
aspects of Lazard's work, starting from the new concept of a
powerful pro-$p$ group. In their setting, uniformly powerful groups
take over the central role which was originally played by saturable
pro-$p$ groups; see \cite{DDMS}. It is easily seen that uniformly
powerful pro-$p$ groups form a subclass of the class of saturable
pro-$p$ groups. Klopsch pointed out that the Sylow pro-$p$ subgroups
of many classical groups are saturable, but typically fail to be
uniformly powerful. Generalising results of Ilani, he established a
one-to-one correspondence between closed subgroups of a saturable
pro-$p$ group and Lie sublattices of the associated $\Z_p$-Lie
lattice; see~\cite{Kl05a}. This somewhat re-established Lazard's
`groupes $p$-saturables' as key players in the Lie theory of
$p$-adic analytic groups and led to applications, for instance, in
the subject of subgroup growth; cf.~\cite{Kl05b}. More recently,
Gonz\'alez-S\'anchez~\cite{GS} gave a characterisation of saturable
pro-$p$ groups, which is more suitable for applications in group
theory than Lazard's original definition: a pro-$p$ group is
saturable if and only if it is finitely generated, torsion-free and
admits a potent filtration; see Section~\ref{sec2} for details. In
general it is a delicate issue to decide whether a given
torsion-free $p$-adic analytic pro-$p$ group is saturable or not.
Using Gonz\'alez-S\'anchez' characterisation, we prove

\begin{teoABC} \label{teo_A} Every torsion-free $p$-adic analytic
  pro-$p$ group of dimension less than $p$ is saturable. On the other
  hand, there exists a torsion-free $p$-adic analytic pro-$p$ group of
  dimension $p$ which is not saturable.
\end{teoABC}

In the present paper Theorem~\ref{teo_A} leads us to several related
results, applications and open questions. For instance, it is the
basic ingredient for an analogue of Lazard's classical correspondence
between finite $p$-groups and finite nilpotent Lie rings of $p$-power
order.

\subsection{Isomorphism of Categories}

As indicated above, a most remarkable and useful property of a
saturable pro-$p$ group is that its underlying set carries naturally
the `dual' structure of a $\Z_p$-Lie lattice. Theorem~\ref{teo_A} and
a similar result for Lie lattices (see Theorem~\ref{teo_A_Lie}) lead
to the following specific instance of Lazard's general Lie
correspondence between saturable pro-$p$ groups and saturable
$\Z_p$-Lie lattices; cf.~\cite[IV~(3.2.6)]{La65} and \cite{Kl05a}. Let
$\G_{<p}$ denote the category of torsion-free $p$-adic analytic
pro-$p$ groups of dimension less than $p$, and let $\L_{<p}$ denote
the category of residually-nilpotent $\Z_p$-Lie lattices of dimension
less than $p$. As we recall in Section~\ref{sec2}, the Hausdorff
series and its inverse give rise to functors $\Grp : \L_{<p}
\rightarrow \G_{<p}$ and $\Lie : \G_{<p} \rightarrow \L_{<p}$.  These
functors can be understood by means of the familiar exponential and
logarithm series; cf.~\cite[III~\S7]{Bo} and \cite[IV~(3.2)]{La65}.

\begin{teoABC}\label{teo_B}
  The functors $\Grp$ and $\Lie$ are mutually inverse isomorphisms of
  categories between $\L_{<p}$ and $\G_{<p}$.
  \begin{enumerate}
  \item The functor $\Grp$ sends $\Z_p$-Lie sublattices to closed
    subgroups, and Lie ideals to closed normal subgroups. Conversely,
    the functor $\Lie$ sends closed subgroups to $\Z_p$-Lie
    sublattices, and closed normal subgroups to Lie ideals.
  \item The functors $\Grp$ and $\Lie$ restrict to isomorphisms
    between the subcategories of nilpotent (respectively soluble)
    objects on both sides. Moreover, they preserve the nilpotency
    class (respectively derived length), when the groups and Lie
    lattices involved are nilpotent (respectively soluble).
  \end{enumerate}
\end{teoABC}

In Section~\ref{sec:classification} we illustrate the usefulness of
this theorem by giving a concrete application: based on an analysis of
conjugacy classes in $\textup{SL}_2(\Z_p)$, we obtain an effective
classification of $3$-dimensional soluble torsion-free $p$-adic
analytic pro-$p$ groups for $p > 3$; see
Theorem~\ref{thm:complete_list}. Our procedure is comparable to the
use of Lazard's correspondence in classifying finite-$p$ groups of
small order, e.g.\ see~\cite{OBVa05}. Implicitly we obtain a complete
list of $3$-dimensional soluble $\Z_p$-Lie lattices; this has
already been used in \cite{KlVo07}. We remark that among the
$3$-dimensional soluble torsion-free $p$-adic analytic pro-$p$ groups
there are several which are not powerful. They complement the
insoluble examples of saturable but not uniformly powerful pro-$p$
groups given in \cite{Kl05a}. For completeness we also discuss briefly
insoluble torsion-free $p$-adic analytic pro-$p$ groups of dimension
$3$.

\subsection{Correspondence between subgroups and Lie sublattices}

In view of Theorem~\ref{teo_A}, the borderline groups of dimension $p$
are of special interest. We prove that saturable pro-$p$ groups of
dimension at most $p$ are in fact saturable in a `strong' sense.

\begin{proABC} \label{pro_C}
  Every saturable pro-$p$ group $G$ of dimension at most $p$ satisfies
  $\gamma_p(G) \subseteq \Phi(G)^p$.
\end{proABC}

We remark that every torsion-free finitely generated pro-$p$ group $G$
satisfying $\gamma_p(G) \subseteq \Phi(G)^p$ is saturable; see
Theorem~\ref{equivalencia}. Moreover, Gonz\'alez-S\'anchez has given
an example of a saturable pro-$p$ group $G$ of dimension $p+1$ such
that $\gamma_p(G) \not \subseteq \Phi(G)^p$; see \cite[Example after
Corollary~3.6]{GS}. Proposition~\ref{pro_C} shows that no such
examples exist in dimensions up to $p$.

\begin{proABC} \label{pro_D}
  Every $p$-adic analytic pro-$p$ group of dimension $p$ which embeds
  into a saturable pro-$p$ group is itself saturable.
\end{proABC}

As an immediate consequence of this inconspicuous proposition we
obtain a new and conceptually satisfying proof of

\begin{teoABC}\label{teo_E}
  Let $G$ be a saturable pro-$p$ group and let $L(G)$ denote the
  associated saturable $\Z_p$-Lie lattice. Let $K, H \subseteq G$ be
  closed subsets, and denote them by $L(K)$, $L(H)$ when regarded as
  subsets of $L(G)$.
  \begin{enumerate}
  \item Suppose that $H$ is a subgroup of $G$ such that every
    $2$-generated subgroup of $H$ has dimension at most $p$. Then
    $L(H)$ is a Lie sublattice of $L(G)$. Moreover, if $K$ is a normal
    subgroup of $H$, then $L(K)$ is a Lie ideal of $L(H)$.
  \item Suppose that $L(H)$ is a Lie sublattice of $L(G)$ such that
    every $2$-generated Lie sublattice of $L(H)$ has dimension at most
    $p$. Then $H$ is a subgroup of $G$. Moreover, if $L(K)$ is a Lie
    ideal of $L(H)$, then $K$ is a normal subgroup of $H$.
  \end{enumerate}
\end{teoABC}

This Lie correspondence between subgroups of a saturable pro-$p$ group
and Lie sublattices of the associated Lie lattice was originally
discovered by Ilani~\cite{Il95} in the context of uniform pro-$p$
groups, and subsequently extended by Klopsch~\cite{Kl05a} to cover the
more general case. Ilani and Klopsch both missed the central insight
which we recorded as Proposition~\ref{pro_D} and had to rely on a less
transparent line of argument.

\subsection{Open questions}

Our results raise several open questions which we record together with
some partial results. In the present paper we work within the class of
all $p$-adic analytic pro-$p$ groups. It would be interesting to
develop saturability criteria, which are tailor-made for more specific
families of groups. Foremost we have in mind the class of $p$-adic
analytic just-infinite pro-$p$ groups, which was systematically
studied in \cite{KlLePl97}. Indeed, the original starting point for
our present work was a result of Klopsch in this direction
\cite[Theorem~1.3]{Kl05a}: every insoluble maximal $p$-adic analytic
just-infinite pro-$p$ group of dimension less than $p-1$ is
saturable. In Section~\ref{sec:just-infinite} we consider

\begin{que} \label{que:just-infinite} What group theoretic properties
  ensure that a $p$-adic analytic just-infinite pro-$p$ group is
  torsion-free or saturable?
\end{que}

In order to formulate a first necessary condition we recall that to
every $p$-adic analytic group $G$ one associates a $\Q_p$-Lie algebra
$\mathcal{L}(G)$; see Section~\ref{sec2}. It is known that the
$\Q_p$-Lie algebra attached to an insoluble $p$-adic analytic
just-infinite pro-$p$ group $G$ is semisimple; in fact, there exists
$k \in \N_0$ such that $\mathcal{L}(G)$ is a direct sum of $p^k$
copies of a simple $\Q_p$-Lie algebra.

\begin{proABC} \label{pro_F} The $\Q_p$-Lie algebra associated to a
  saturable insoluble just-infinite pro-$p$ group is simple. Every
  saturable insoluble just-infinite pro-$p$ group is hereditarily
  just-infinite.

  In contrast, for any simple $\Q_p$-Lie algebra $\mathcal{L}$ and for
  any $p$-power $p^k$, $k \in \N_0$, there exists a torsion-free
  $p$-adic analytic just-infinite pro-$p$ group whose associated
  $\Q_p$-Lie algebra is isomorphic to the direct sum of $p^k$ copies
  of $\mathcal{L}$.
\end{proABC}

In particular, the proposition yields just-infinite examples of
non-saturable torsion-free $p$-adic analytic pro-$p$ groups of
dimension $3p$. This indicates natural limitations of a possible
extension of Theorem~\ref{teo_A}.

A trivial consequence of Proposition~\ref{pro_D} is that saturable
pro-$p$ groups of dimension at most $p$ have the property that every
one of their open subgroups is again saturable. It is a challenging
problem to find other natural families of pro-$p$ groups which are
`hereditarily saturable'. This motivates

\begin{que}\label{que:all_open_saturable}
  Are all open subgroups of a saturable just-infinite pro-$p$ group
  saturable?
\end{que}

In view of~\cite[Theorem~1.1]{Kl05a} we may specialise our question
further to special linear groups over the ring of integers of a
$p$-adic field as follows. Let $\mathcal{O}$ be the ring of integers
of a finite extension of $\Q_p$ with ramification index $e$, and let
$n \in \N$ with $n e \leq p-2$. Are all open subgroups of a Sylow
pro-$p$ subgroup of $\textup{SL}_n(\mathcal{O})$ saturable?

\smallskip

In a different direction, it would be interesting to explore further
the borderline cases where Lie theoretic methods break down. In light
of Theorem~\ref{teo_A} it is natural to ask
\begin{que} \label{que:dim_p} What are the torsion-free $p$-adic
  analytic pro-$p$ groups $G$ of dimension $p$ which are not
  saturable, i.e.\ which fail to satisfy the condition $\gamma_p(G)
  \subseteq \Phi(G)^p$? Can we classify them? Similarly, what are the
  residually-nilpotent $\Z_p$-Lie lattices of dimension $p$ which are
  not saturable?
\end{que}

Based on the examples that we have at present, it is conceivable that
one can go beyond a Lie theory founded on Lazard's notion of
saturability. Pink has interesting results which support this
idea; cf.~\cite{Pi93}.

\begin{que} \label{que:Lie_cor}
  Is there a (natural) correspondence between residually-nilpotent
  $\Z_p$-Lie lattices of dimension $p$ and torsion-free $p$-adic
  analytic pro-$p$ groups of dimension $p$?
\end{que}

\smallskip

\noindent
\textit{Organisation.} In Section~\ref{sec2} we recall the
characterisations of saturable pro-$p$ groups and saturable $\Z_p$-Lie
lattices in terms of potent filtrations. Subsequently, we summarise
Lazard's correspondence between saturable pro-$p$ groups and saturable
$\Z_p$-Lie lattices. In Section~\ref{sec:iso} we prove an auxiliary
result concerning isolated subgroups of saturable pro-$p$ groups. In
Section~\ref{sec:dim_p} we prove Theorem~\ref{teo_A},
Proposition~\ref{pro_C} and Proposition~\ref{pro_D}. As immediate
consequences we obtain Theorems~\ref{teo_E} and \ref{teo_B}. In
Section~\ref{sec:structure} we provide in the context of $\Z_p$-Lie
lattices analogues of the classical Levi splitting and the
decomposition of semisimple Lie algebras into simple summands. From
this we prove the first part of Proposition~\ref{pro_F}. In
Section~\ref{sec:just-infinite} we look at $p$-adic analytic
just-infinite pro-$p$ groups and complete the proof of
Proposition~\ref{pro_F}. In Section~\ref{sec:classification} we
describe torsion-free $p$-adic analytic pro-$p$ groups of dimension at
most $3$ by a procedure which is comparable to the use of Lazard's
correspondence in classifying finite-$p$ groups of small order. Our
main result is Theorem~\ref{thm:complete_list} which provides an
effective classification of $3$-dimensional soluble torsion-free
$p$-adic analytic pro-$p$ groups for $p > 3$.

\smallskip

\noindent
\textit{Notation.} Throughout, $p$ denotes an odd prime. Let $G$ be a
group. If $H \subseteq G$ and $n \in \N$, then $H^n := \langle h^n
\mid h \in H \rangle$. The Frattini subgroup of $G$ is denoted by
$\Phi(G)$. In particular, if $G$ is a pro-$p$ group, then $\Phi(G) =
G^p [G,G]$. As customary, the terms of the lower central series of $G$
are denoted by $\gamma_i(G)$, $i \in \N$.

As a rule we do not distinguish notationally between the group
commutator $[x,y]$ of group elements $x,y$ and the Lie bracket $[x,y]$
of Lie lattice elements $x,y$. Both, group commutators and Lie
brackets are left-normed. For $k \in \mathbb{N}$, the abbreviation
$[N,_k M]$ stands for $[N,M,\ldots ,M]$ with $M$ occurring $k$ times.


\section{Potent filtrations and Lie correspondence} \label{sec2}

In the present section we recall Gonz\'alez-S\'anchez'
characterisations of saturable pro-$p$ groups and saturable Lie
lattices, based on the concept of potent filtrations. Then we
summarise Lazard's correspondence between saturable pro-$p$ groups and
saturable $\Z_p$-Lie lattices.

Let $G$ be a pro-$p$ group, and let $N$ be a closed normal subgroup of
$G$.  A \emph{potent filtration} of $N$ in $G$ is a descending series
$(N_i)_{i \in \N}$ of closed normal subgroups of $G$ such that (i)
$N_1 = N$, (ii) $\bigcap_{i \in \N} N_i = 1$, (iii) $[N_i,G] \subseteq
N_{i+1}$ and $[N_i,_{p-1} G] \subseteq N_{i+1}^p$ for all $i \in \N$.
We say that $N$ is \emph{PF-embedded} in $G$ if there exists a potent
filtration of $N$ in $G$. The group $G$ is a \emph{PF-group}, if $G$
is PF-embedded in itself.

We recall basic properties of PF-embedded subgroups, which follow
essentially from the Hall-Petrescu collection formula~\cite[\S9]{Hu};
see~\cite[Proposition~3.2 and Theorem~3.4]{FGJ} and
\cite[Proposition~2.2]{GS}.

\begin{lem}[Fern\'andez-Alcober, Gonz\'alez-S\'anchez,
  Jaikin-Zapirain] \label{propiedades}
  Let $G$ be a pro-$p$ group, and let $N, M$ be PF-embedded subgroups
  of $G$. Then
  \begin{enumerate}
  \item $NM$, $N^p$ and $[N,_kG]$ are PF-embedded in $G$ for all $k
    \in \N$;
  \item $[N^p,G] = [N,G]^p$;
  \item $N^p = \{ x^p \mid x \in N \}$;
  \item if $G$ is torsion-free and $x \in G$ with $x^p \in N^p$, then
    $x \in N$; moreover, if $x,y \in N$ such that $x^p = y^p$, then $x
    = y$.
  \end{enumerate}
\end{lem}

We also state Gonz\'alez-S\'anchez' characterisation of saturable
pro-$p$ groups; see~\cite[Theorem~3.4 and Corollary~5.4]{GS}.

\begin{teo}[Gonz\'alez-S\'anchez]\label{equivalencia}
  Let $G$ be a torsion-free finitely generated pro-$p$ group. Then $G$
  is saturable if and only if $G$ (or equivalently $G /
  \Phi(G)^p$) is a PF-group.

  In particular, if $\gamma_p(G) \subseteq \Phi(G)^p$, then $G$ is
  saturable.
\end{teo}

If $G$ is a torsion-free finitely generated pro-$p$ group satisfying
$\gamma_p(G) \subseteq \Phi(G)^p$, then it is particularly easy to
write down a potent filtration $G_i$, $i \in \N$: one simply takes the
lower $p$-series
$$
G_1 := G, \quad \text{and} \quad G_i := [G_{i-1},G] G_{i-1}^p \text{
  for $i \geq 2$.}
$$
Indeed, $[G_1,_{p-1} G] = \gamma_p(G) \subseteq \Phi(G)^p = G_2^p$
and, using induction and Lemma~\ref{propiedades}, we see that for $i
\geq 2$,
\begin{align*}
  [G_i,_{p-1} G] & = [[G_{i-1},G] G_{i-1}^p,_{p-1} G]
  \subseteq [G_{i-1},_{p-1}G,G] [G_{i-1}^p,_{p-1} G] \\
  & \subseteq [G_i^p,G] [G_{i-1},_{p-1} G]^p \subseteq [G_i,G]^p
  (G_i^p)^p \subseteq G_{i+1}^p.
\end{align*}

Similarly as for groups, one can study saturable $\Z_p$-Lie lattices
and Lie lattices admitting a potent filtration by ideals; for details
see \cite[I~(2)]{La65}, summarised in \cite[\S2]{Kl05a}, and
\cite[\S4]{GS}. The analogue of Theorem~\ref{equivalencia} states that
a $\Z_p$-Lie lattice $L$ is saturable if and only if $L$ (or
equivalently $L/(p[L,L]+p^2L)$) admits a potent filtration.

Lazard's correspondence between saturable pro-$p$ groups and saturable
$\Z_p$-Lie lattices can be obtained via exponential and logarithm maps
inside certain Hopf algebras. Alternatively, the passage from Lie
lattice to group structure can be described in terms of the Hausdorff
series (sometimes called Baker-Campbell-Hausdorff series,
cf.~\cite[Historical Note,~V]{Bo}) and its inverse.

The Hausdorff series is defined as $\Phi(X,Y) := \log((\exp X)(\exp
Y)) \in \Q\langle\!\langle X,Y \rangle\!\rangle$. For our purposes
one writes $\Phi(X,Y)$ as a series of Lie words in independent
variables $X,Y$. For instance, modulo terms of weight $4$ or higher,
we have
$$
\Phi(X,Y) \equiv X + Y + \frac{1}{2} [X,Y] + \frac{1}{12} ( [X,Y,Y] -
[X,Y,X] );
$$
cf.~\cite[II~\S8]{Bo}. Given a saturable $\Z_p$-Lie lattice $L$, the
saturable pro-$p$ group corresponding to $L$ has underlying set $L$
and the group product of $x,y \in L$ is given by $x y :=
\Phi(x,y)$. Conversely, given a saturable pro-$p$ group $G$, one
defines on the set $G$ the structure of a saturable $\Z_p$-Lie lattice
$L(G)$ as follows. For $\lambda \in \Z_p$ and $x,y \in G$ one sets
$$
    \lambda x := x^\lambda, \quad
    x + y := \lim_{n \rightarrow \infty} (x^{p^n}
    y^{p^n})^{p^{-n}}, \quad
    [x,y]_\textup{Lie} := \lim_{n \rightarrow \infty}
    [x^{p^n},y^{p^n}]^{p^{-2n}}.
$$
See~\cite[IV~(3.2.6)]{La65} or \cite[\S2]{Kl05a}; the latter contains
a more detailed summary. The $\Q_p$-Lie algebra associated to a
compact $p$-adic analytic group $G$ is define as $\mathcal{L}(G) :=
\Q_p \otimes_{\Z_p} L(U)$, where $U$ is an open saturable pro-$p$
subgroup of $G$. Clearly, $\mathcal{L}(G)$ only depends on the
commensurability class of $G$.


\section{Isolators in saturable pro-$p$ groups} \label{sec:iso}

An effective tool, used in \cite{Kl05a}, consists in embedding a given
subgroup of a saturated pro-$p$ group into its saturated closure; also
cf.\ \cite[IV~(3.2.6)]{La65}. We now explain this procedure using
potent filtrations.

Let $G$ be a closed subgroup of a pro-$p$ group $S$. The
\emph{isolator} of $G$ in $S$ is defined as the closed subgroup
$$
\iso_S(G) := \langle g \in S \mid \exists k \in \N : g^{p^k} \in G \rangle.
$$
Now suppose that $S$ is saturable. Then Theorem~\ref{equivalencia}
shows that $S$ admits a potent filtration $S_i$, $i \in \N$. Note that
$S_{i+1}^p = S_{i+1}^{\{p\}}$ for all $i \in \N$, by
Lemma~\ref{propiedades}, and thus
$$
[(\iso_S(G) \cap S_i),_{p-1} \iso_S(G)] \subseteq \iso_S(G) \cap
S_{i+1}^{\{p\}} = (\iso_S(G) \cap S_{i+1})^{\{p\}}.
$$
This implies that $\iso_S(G) \cap S_i$, $i \in \N$, is a potent
filtration of $\iso_S(G)$. By Theorem~\ref{equivalencia}, the
group $\iso_S(G)$ is saturable.

If $G$ is a normal subgroup of $S$, the isolator $\iso_S(G)$ is also a
normal subgroup of $S$, in particular
$$
[(\iso_S(G) \cap S_i),_{p-1} S] \subseteq \iso_S(G) \cap
S_{i+1}^{\{p\}} = (\iso_S(G) \cap S_{i+1})^{\{p\}}.
$$
Hence $\iso_S(G) \cap S_i$, $i\in N$, is a potent filtration of $S$
starting at $\iso_S(G)$, and thus $\iso_S(G)$ is PF-embedded in $S$.

At this stage it is advantageous to briefly consider the analogous
situation for $\Z_p$-Lie lattices. Let $L$ be a Lie sublattice of a
$\Z_p$-Lie lattice $S$. Then the isolator of $L$ in $S$ is the Lie
sublattice $\iso_S(L) = S \cap (\Q_p \otimes L)$. We observe that
the concept is a priori easier to handle than for groups. For
instance, it is immediate that $\lvert \iso_S(L) : L \rvert <
\infty$. Similar arguments as the ones given above in the context of
groups lead to

\begin{pro} \label{isolatorLie} Let $L$ be a Lie sublattice of a
  saturable $\Z_p$-Lie lattice $S$. Then the isolator $\iso_S(L)$ is
  saturable and $\lvert \iso_S(L) : L \rvert < \infty$. In particular,
  $L$ and $\iso_S(L)$ have the same dimension.

  Furthermore, if $L$ is a Lie ideal of $S$, then $\iso_S(L)$ is
  PF-embedded in $S$.
\end{pro}

Our aim is to prove the corresponding result for groups, and we return
to the situation where $G$ is a subgroup of a saturable pro-$p$ group
$S$. Recall that the group $S$ and its associated $\Z_p$-Lie lattice
$L(S)$ are defined on the same underlying set and that taking $p$th
powers in $S$ has the same effect as multiplying by $p$ in $L(S)$;
cf.\ Section~\ref{sec2}. The group $G$ contains an open subgroup
$H$ which is saturable, and clearly $\iso_S(G) = \iso_S(H)$. By
Lazard's correspondence $L(H)$, i.e.\ $H$ considered as a subset of
$L(G)$, forms a Lie sublattice of $L(S)$. By
Proposition~\ref{isolatorLie}, the Lie isolator $\iso_{L(S)}(L(H))$ is
saturable and corresponds to a subgroup of $S$. Indeed, this latter
subgroup is easily identified as the group isolator $\iso_S(H)$. We
summarise our findings in

\begin{pro}\label{isolator}
  Let $G$ be a closed subgroup of a saturable pro-$p$ group $S$. Then
  $\iso_S(G) = \{ g \in S \mid \exists k \in \N : g^{p^k} \in G \}$.
  Moreover, $\iso_S(G)$ is saturable and $\lvert \iso_S(G):G \rvert <
  \infty$. In particular, $G$ and $\iso_S(G)$ have the same dimension.

  Furthermore, if $G$ is a normal subgroup of $S$, then $\iso_S(G)$ is
  PF-embedded in $S$.
\end{pro}


\section{Groups of dimension at most $p$} \label{sec:dim_p} In the
present section we prove Theorem~\ref{teo_A}, Proposition~\ref{pro_C}
and Proposition~\ref{pro_D}. As immediate consequences we obtain
Theorems~\ref{teo_E} and \ref{teo_B}.

\begin{pro}\label{p-1}
  Let $G$ be a torsion-free $p$-adic analytic pro-$p$ group of
  dimension less than $p$. Then $\gamma_p(G) \subseteq \Phi(G)^p$; in
  particular, $G$ is a PF-group.
\end{pro}

\begin{proof}
  If $p = 2$, then $\Z_2$ is the only $1$-dimensional torsion-free
  pro-$p$ group and the proposition is obviously true. Let us suppose
  that $p\geq 3$.  By Theorem~\ref{equivalencia}, the assertion
  $\gamma_p(G) \subseteq \Phi(G)^p$ which we are about to prove will
  imply that $G$ is a PF-group. Clearly, we may assume that $G$ is not
  the trivial group.  Since $G$ is $p$-adic analytic, there exists a
  proper open normal subgroup $N$ of $G$ which is saturable and hence
  a PF-group.  The chief factors of the finite $p$-group $G/N$ are
  cyclic of order $p$.  Arguing by induction on the composition length
  of $G/N$, we may assume that $|G:N| = p$.

  Since $N$ is saturable, we have $|N : N^p| = p^{\dim(N)} =
  p^{\dim(G)} \leq p^{p-1}$; see~\cite[III~(3.1)]{La65}. This shows that
  $|G : N^p| \leq p^p$, and consequently $\gamma_p(G) \subseteq N^p$.

  \smallskip

  \noindent \textit{Case 1:} $\gamma_{p-1}(G) \subseteq N^p$. Then, a
  fortiori, $\gamma_{p-1}(G) \leq G^p$ and the Hall-Petrescu
  formula~\cite[\S9]{Hu} shows that, modulo $[G,G]^p$,
  $$
  \gamma_p(G) = [\gamma_{p-1}(G),G] \subseteq [G^p,G] \subseteq
  \gamma_{p+1}(G);
  $$
  cf.~\cite[proof of Theorem~3.1]{GJ}. Thus $\gamma_p(G) \subseteq
  [G,G]^p \subseteq \Phi(G)^p$, as wanted.

  \smallskip

  \noindent \textit{Case 2:} $\gamma_{p-1}(G) \not \subseteq N^p$.
  Then $G/N^p$ is a finite $p$-group of maximal class: $G/N^p$ has
  order $p^p$ and nilpotency class $p-1$. Put $N_1 := N$ and $N_i :=
  \gamma_i(G)N^p$ for $i \in \{2, \ldots, p-1 \}$. Then
  $$
  G \supset N = N_1 \supset N_2 \supset \ldots \supset N_{p-1}
  \supset N^p
  $$
  is a central descending series of closed normal subgroups such
  that each term has index $p$ in its predecessor. Moreover, every
  normal subgroup of $G$ which is strictly contained in $N$ must be
  contained in $N_2 = \gamma_2(G) N^p$.

  First suppose that $\gamma_p(G) \subsetneqq N^p$. The set of
  elements of order $p$ in the regular $p$-group $N / \gamma_p(G)$
  forms a subgroup. Hence $M := \{ x \in N \mid x^p \in \gamma_p(G)
  \}$ forms a subgroup of $G$ which, by assumption, is strictly
  contained in $N$. By our observation, $M$ is thus contained in
  $N_2$. Furthermore, according to Lemma~\ref{propiedades}~(3), we
  have $N^p = \{ x^p \mid x \in N \}$, and hence $M^p = \gamma_p(G)$.
  This gives $\gamma_p(G) = M^p \subseteq N_2^p \subseteq \Phi (G)^p$,
  as wanted.

  We claim that the remaining case $\gamma_p(G) = N^p$ actually does
  not occur. For a contradiction, suppose that $\gamma_p(G) = N^p$. We
  define recursively $N_i := N_{i-p+1}^p$ for $i \in \N$ with $i \geq
  p$. Note that $N_i = \gamma_i(G)$ for $i \in \{2,\ldots, p \}$.

  \smallskip

  \noindent \textit{Subclaim:} $N_i$ is PF-embedded in $N$ for all
  $i \in \N$.

  \noindent \textit{Subproof.} This is true for $i = 1$ as $N = N_1$
  is a PF-group, and by Lemma~\ref{propiedades}~(1) it now suffices to
  prove the assertion for $i \in \{2, \ldots, p-1 \}$. Let $i$ be in
  this range. It follows from Lemma~\ref{propiedades}~(1),(2) that
  $[N^p,_{i-1} N] \subseteq \gamma_i(N)^p$. This gives
  $$
  [N_i,_{p-1} N] = [\gamma_i(G),_{p-1} N] \subseteq
  [\gamma_p(G),_{i-1} N] = [N^p,_{i-1} N] \subseteq \gamma_i(N)^p.
  $$
  Now, $\gamma_i(N)$ is PF-embedded in $N$ by
  Proposition~\ref{propiedades}~(1), hence we can complete the
  inclusion $N_i = \gamma_i(G) \supseteq \gamma_i(N)$ to a potent
  filtration of $N_i$ in $N$. This proves the subclaim.

  \smallskip

  \noindent \textit{Subclaim:} $|N_i : N_{i+1}| = p$ for all $i \in
  \N$.

  \noindent \textit{Subproof.}  This is clear for $i \in \N$ with $i
  \leq p - 1$. Now let $i \geq p$ and argue by induction. By the first
  subclaim $N_j$ is a PF-group for every $j \in \N$; thus
  Lemma~\ref{propiedades}~(3) shows that $N_{j + p - 1} = N_j^p = \{
  x^p \mid x \in N_j \}$. From \cite[Lemma~A.3]{Kl05a} it follows by
  induction that $|N_i : N_{i+1}| = |N_{i-p+1} : N_{i-p+2}| = p$.

  \smallskip

  \noindent \textit{Subclaim:} $N_i = \gamma_i(G)$ for all $i \in
  \N$ with $i \geq 2$.

  \noindent \textit{Subproof.} Let $i \in \N$ with $i \geq 2$, and
  argue by induction. For $i \leq p$ the assertion holds, as noted
  above. Suppose now that $i > p$. By induction we have $N_i =
  N_{i-p+1}^p = \gamma_{i-p+1}(G)^p = [\gamma_{i-p}(G),G]^p =
  [N_{i-p},G]^p$ and similarly $\gamma_i(G) = [\gamma_{i-1}(G),G] =
  [N_{i-1},G] = [N_{i-p}^p,G]$. The Hall-Petrescu
  formula~\cite[\S9]{Hu} yields $[N_{i-p},G]^p \equiv [N_{i-p}^p,G]$
  modulo $[G,_pN_{i-p}]$ (cf.~ \cite[Theorem~2.4]{FGJ}), thus
  $$
  N_i \equiv \gamma_i(G) \pmod{[G,_p N_{i-p}]}.
  $$

  By induction, $[G,_p N_{i-p}] = [G,_p \gamma_{i-p}(G)] \subseteq
  \gamma_i(G)$. Again by induction and since $N_{i-p+1}$ is
  PF-embedded in $N$, we also obtain
  \begin{align*}
    [G,_p N_{i-p}] & = [G, \gamma_{i-p}(G),_{p-1} N_{i-p}] \subseteq
    [\gamma_{i-p+1}(G),_{p-1} N] \\
    & \subseteq [N_{i-p+1},_{p-1} N] \subseteq N_{i-p+2}^p = N_{i+1}
    \subseteq N_i.
  \end{align*}
  Therefore it follows that $N_i = \gamma_i(G)$, as claimed.

  \smallskip

  The last two subclaims show that $G$ is the unique pro-$p$ group of
  maximal class; cf.~\cite[Section~7.4]{LeMc02}. But in this case $G$
  has elements of order $p$, a contradiction.
\end{proof}

The following example indicates that the condition on the dimension in
Proposition~\ref{p-1} is sharp.

\begin{eje}\label{p}
  We display a torsion-free $p$-adic analytic pro-$p$ group of
  dimension $p$ which is not a PF-group.

  Let $M = \langle x_1,\ldots ,x_{p-1}\rangle \cong \Z_p^{p-1}$ be a
  free abelian pro-$p$ group of rank $p-1$, and let $A = \langle
  \alpha \rangle \cong \Z_p$. Consider the semidirect product $G := A
  \ltimes M$, with respect to the action of $A$ on $M$ defined by
  $$
  x_i^\alpha =
  \begin{cases}
    x_i x_{i+1} & \text{if $1 \leq i \leq p-2$,} \\
    x_{p-1} x_1^p & \text{if $i = p-1$.}
  \end{cases}
  $$
  Clearly, $G$ is a torsion-free $p$-adic analytic pro-$p$ group.
  Moreover, it is easily seen that $[M,_{p-1} G] = M^p$. For a
  contradiction, suppose that $G$ is a PF-group and let $( G_i )_{i
    \in \N}$ be a potent filtration of $G$. Then $M^p \subseteq
  [G_1,_{p-1} G] \subseteq G_2^p$ implies that $M \subseteq G_2$.
  Inductively, it follows that $M \subseteq G_i$ for all $i \in \N$.
  This contradicts the fact that $\bigcap G_i = 1$. \hfill $\Diamond$
\end{eje}

Theorem~\ref{teo_A} follows directly from Proposition~\ref{p-1} and
Example~\ref{p}.

\begin{cor}\label{p-1-PF}
  Let $G$ be a saturable pro-$p$ group, and let $N$ be a normal
  subgroup of $G$ of dimension less than $p$. Then $N$ is PF-embedded
  in $G$.
\end{cor}

\begin{proof}
  By Theorem~\ref{equivalencia}, $G$ is a PF-group, and by
  Proposition~\ref{p-1}, $N$ is a PF-group and hence saturable.
  Let $(G_i)_{i\in \N}$ be a potent filtration of $G$. Define $N_i :=
  N \cap G_i$ for $i \in \N$. Since $N$ is saturable of dimension
  less than $p$, we have $|N:N^p| \leq p^{p-1}$; see
  \cite[III~(3.1)]{La65}.  Therefore $[N,_{p-1} G] \subseteq N^p$, and
  hence $[N_i,_{p-1} G] \subseteq N^p \cap G_{i+1}^p$. Applying
  Lemma~\ref{propiedades}~(3),(4), we obtain $N^p \cap G_{i+1}^p = (N
  \cap G_{i+1})^p = N_{i+1}^p$ for all $i \in \N$. Thus $(N_i)_{i \in
    \N}$ is a potent filtration of $N$ in $G$.
\end{proof}

Proposition~\ref{p-1} implies that every saturable pro-$p$ group $G$
of dimension less than $p$ satisfies $\gamma_p(G) \subseteq
\Phi(G)^p$.  In order to complete the proof of
Proposition~\ref{pro_C}, it suffices to show

\begin{pro}
  Let $G$ be a saturable pro-$p$ group of dimension $p$. Then
  $\gamma_p(G)\subseteq \Phi(G)^p$.
\end{pro}

\begin{proof}
  Consider the finite $p$-group $H := G/\Phi(G)^p \gamma_{p+1}(G)$. By
  definition this is a PF-group of nilpotency class at most $p$, and
  our task is to prove that, in fact, $\gamma_p(H) = 1$

  Note that $|H : H^p| \leq |G : G^p| = p^{\dim(G)} = p^p$. This shows
  that $\gamma_p(H) \subseteq H^p$. First suppose that $|H:\Phi(H)|
  \geq p^3$. As $\Phi(H) = \gamma_2(H) H^p$, we deduce that
  $\gamma_{p-1}(H) \subseteq H^p$, and Lemma~\ref{propiedades}~(2)
  implies
  $$
  \gamma_p(H) = [\gamma_{p-1}(H),H] \subseteq [H^p,H] = [H,H]^p
  \subseteq \Phi(H)^p = 1,
  $$
  as wanted.

  Now suppose that $|H:\Phi(H)| = p^2$, i.e.\ that $H$ is
  $2$-generated. For a contradiction, we assume that $\gamma_p(H) \ne
  1$. Choose $z \in \gamma_p(H) \setminus \{1\}$. Since $\gamma_p(H)
  \subseteq H^p$, Lemma~\ref{propiedades}~(3) shows that
  $$
  Y = \{ y \in H \mid y^p = z \} \ne \varnothing.
  $$
  Since $Y \cap \Phi(H) = \varnothing$, every $y \in Y$ forms part of
  a generating pair of $H$. Lemma~\ref{propiedades}~(2) shows that
  $[H^p,_{p-1} H] \subseteq [H,H]^p \subseteq \Phi(H)^p = 1$, hence
  $$
  [\Phi(H),_{p-1} H] \subseteq [H^p,_{p-1}H] [[H,H],_{p-1} H] =
  \gamma_{p+1}(H) = 1.
  $$ This implies that $\gamma_p(H) = [y,_{p-1} H]$ for every $y \in
  Y$.

  Let $(H_i)_{i \in \N}$ be a potent filtration of $H$. To obtain the
  desired contradiction, we show that $Y \cap H_i \ne \varnothing$ for
  all $i \in \N$. Clearly, the assertion is true for $i = 1$. Now
  consider $i \geq 2$. By induction, we find $y \in Y \cap H_{i-1}$,
  and hence $z \in \gamma_p(H) = [y,_{p-1} H] \subseteq
  [H_{i-1},_{p-1} H] = H_i^p$. Hence Lemma~\ref{propiedades}~(3)
  implies that $Y \cap H_i \ne \varnothing$.
\end{proof}

The Lie correspondence stated as Theorem~\ref{teo_E} connects
subgroups and Lie sublattices of dimension at most $p$. A satisfying
explanation for the underlying phenomenon comes from the next
result, which we stated as Proposition~\ref{pro_D} in the
introduction.

\begin{cor}
  Every $p$-adic analytic pro-$p$ group of dimension $p$ which embeds
  into a saturable pro-$p$ group is itself saturable.
\end{cor}

\begin{proof}
  Let $G$ be a subgroup of a saturable pro-$p$ group $S$ and suppose
  that $\dim(G) \leq p$. By Proposition~\ref{isolator} we may assume
  that $G$ has finite index in $S$ and by induction we can further
  assume that $|S:G| = p$.  Proposition~\ref{pro_C} then shows that
  $\gamma_p(S) \subseteq \Phi(S)^p$. As explained after
  Theorem~\ref{equivalencia}, the lower $p$-series
  $$
  S_1 := S, \quad \text{and} \quad S_i := [S_{i-1},S] S_{i-1}^p \text{
    for $i \geq 2$}
  $$
  is a potent filtration of $S$. Note that $|S:G| = p$ implies $S_2 =
  \Phi(S) \subseteq G$. Hence $G \supset S_2 \supset S_3 \supset \ldots$
  supplies a potent filtration of $G$ and $G$ is saturable.
\end{proof}

The results in the present section, notably Proposition~\ref{p-1} and
Corollary~\ref{p-1-PF}, have counterparts for $\Z_p$-Lie lattices; the
corresponding proofs are very similar to the ones for pro-$p$ groups
-- in fact, they are somewhat easier.

\begin{teo}\label{teo_A_Lie}
  Every residually-nilpotent $\Z_p$-Lie lattice of dimension less than
  $p$ is saturable.

  Furthermore, if $L$ is a saturable $\Z_p$-Lie lattice and $M$ is a Lie
  ideal of $L$ of dimension less than $p$, then $M$ is PF-embedded in
  $L$.
\end{teo}

In view of \cite[\S4]{GS}, we can now combine the last theorem,
Proposition~\ref{p-1} and Corollary~\ref{p-1-PF} to deduce: there is a
correspondence between torsion-free $p$-adic analytic pro-$p$ groups
of dimension less than $p$ and residually-nilpotent $\Z_p$-Lie
lattices of dimension less than $p$ with the additional properties
stated in Theorem~\ref{teo_B}.

For completeness we record the analogue of Example~\ref{p}.
\begin{eje}\label{lie p}
  We display a $\Z_p$-Lie lattice of dimension $p$ which is not
  saturable. Consider the Lie lattice $L = \Z_p x + \Z_p y_1
  + \ldots + \Z_p y_{p-1}$ of dimension $p$ where $[y_i,y_j] = 0$
  for $1 \leq i,j \leq p-1$ and
  $$
  [y_i,x] =
  \begin{cases}
    y_{i+1} & \text{if $i<p-1$,}\\
    p y_1 & \text{if $i=p-1$.}
  \end{cases}
  $$
  It is easy to check that, indeed, $L$ is not saturable.

  Furthermore, it is impossible to associate a pro-$p$ group to the
  Lie lattice $L$ by direct application of the Hausdorff series.  The
  problem in this specific case is that, although all homogeneous
  summands which appear in the formula exist, their sum does not
  necessarily converge (for example, try to compute $\Phi(x,y_1)$).
  \hfill $\Diamond$
\end{eje}

Note that the Lie lattice in Example~\ref{lie p} and the pro-$p$ group
in Example~\ref{p} are constructed in a very similar way. This
observation leads naturally to Questions \ref{que:dim_p} and
\ref{que:Lie_cor} in the introduction.


\section{The structure of saturable groups} \label{sec:structure}

In the present section we apply the notion of isolators to derive for
$\Z_p$-Lie lattices analogues of the classical Levi splitting and the
decomposition of semisimple Lie algebras into simple
summands. Subsequently, we relate these results to saturable pro-$p$
groups and prove the first part of Proposition~\ref{pro_F}.

The \emph{soluble radical} of a $\Z_p$-Lie lattice $L$ is the unique
maximal soluble Lie ideal of $L$. There is a close connection between the
soluble radical of $L$ and the soluble radical of the $\Q_p$-Lie
algebra $\Q_p \otimes L$.

\begin{pro} \label{solvable_radical} Let $L$ be a $\Z_p$-Lie lattice.
  Then the soluble radical $R$ of $L$ is isolated, i.e.\ $R =
  \iso_L(R)$. Moreover, if $L$ is saturable, then $R$ is PF-embedded
  in $L$.
\end{pro}

\begin{proof}
  Clearly, the intersection of $L$ with any soluble Lie ideal of the
  $\Q_p$-Lie algebra $\Q_p \otimes L$ yields a soluble Lie ideal of
  $L$.  Conversely, if $I$ is a soluble Lie ideal of $L$, then $\Q_p
  \otimes I$ constitutes a soluble Lie ideal of $\Q_p \otimes
  L$. Hence the soluble radical of $\Q_p \otimes L$ is equal to $\Q_p
  \otimes R$, and $R = L \cap (\Q_p \otimes R)$ is isolated. The final
  claim follows from Proposition~\ref{isolatorLie}.
\end{proof}

The classical Levi splitting of $\Q_p$-Lie algebras (cf.\ \cite[Part
I, VI.4]{Se65}) yields a corresponding splitting of $\Z_p$-Lie
lattices up to finite index.

\begin{pro}
  Let $L$ be a $\Z_p$-Lie lattice with soluble radical $R$.  Then
  there exists a Lie sublattice $H$ of $L$ with $H \cap R = 0$ such
  that the semidirect sum $H + R$ has finite index in $L$.
\end{pro}

\begin{proof}
  Let $\mathcal{R}$ be the soluble radical of the $\Q_p$-Lie algebra
  $\mathcal{L} := \Q_p \otimes L$. By Levi's classical splitting
  theorem, $\mathcal{L}$ is the semidirect sum of $\mathcal{R}$ and a
  suitable subalgebra $\mathcal{H}$. Write $H := \mathcal{H} \cap L$
  and observe that $R = \mathcal{R} \cap L$ is the soluble radical of
  $L$.  Then clearly $H \cap R = 0$ and, comparing dimensions, we see
  that the Lie sublattice $H + R$ has finite index in $L$.
\end{proof}

\begin{eje}
  We display a powerful $\Z_p$-Lie lattice which is not the semidirect
  sum of a suitable Lie sublattice and the soluble radical. This
  illustrates that in general one cannot hope for a stronger form of
  Levi splitting in $\Z_p$-Lie lattices corresponding to saturable
  pro-$p$ groups.

  Let $k \in \N$ with $k \geq 2$, and let $L$ be the $5$-dimensional
  $\Z_p$-Lie sublattice of $\mathfrak{gl}_3(\Z_p)$ with $\Z_p$-basis
  $$
  \mathbf{x} := p^k
  \begin{pmatrix}
    0 & 1 & 0 \\
    0 & 0 & 0 \\
    1 & 0 & 0
  \end{pmatrix}, \quad \mathbf{y} := p^k
  \begin{pmatrix}
    0 & 0 & 0 \\
    1 & 0 & 0 \\
    0 & 0 & 0
  \end{pmatrix}, \quad \mathbf{h} := p^k
  \begin{pmatrix}
    1 & 0  & 0 \\
    0 & -1 & 0 \\
    0 & 0 & 0
  \end{pmatrix},
  $$
  $$
  \mathbf{a} := p^{2k-1}
  \begin{pmatrix}
    0 & 0 & 0 \\
    0 & 0 & 0 \\
    1 & 0 & 0
  \end{pmatrix}, \quad \mathbf{b} := p^{2k-1}
  \begin{pmatrix}
    0 & 0 & 0 \\
    0 & 0 & 0 \\
    0 & 1 & 0
  \end{pmatrix}.
  $$
  Clearly, $[L,L] \subseteq pL$ and the Lie lattice $L$ is
  powerful. The $\Q_p$-Lie algebra $\Q_p \otimes L$ is, of course, the
  semidirect sum of $\mathfrak{sl}_2(\Q_p)$ and $\Q_p \oplus \Q_p$
  with respect to the natural action. Accordingly the radical of $L$
  is $R = \Z_p \mathbf{a} + \Z_p \mathbf{b}$. Let
  $\widetilde{\mathbf{h}} \equiv \mathbf{h}$ and
  $\widetilde{\mathbf{x}} \equiv \mathbf{x}$ modulo $R$.  Then
  $[\widetilde{\mathbf{h}}, \widetilde{\mathbf{x}}] \equiv
  [\mathbf{h},\mathbf{x}]$ modulo $p^kR$. Moreover,
  $[\mathbf{h},\mathbf{x}] - 2p^k \widetilde{\mathbf{x}} \equiv 0$
  modulo $R$, but $[\mathbf{h},\mathbf{x}] - 2p^k
  \widetilde{\mathbf{x}} \not \equiv 0$ modulo $p^k R$. This shows
  that the Lie sublattice generated by $\widetilde{\mathbf{x}}$ and
  $\widetilde{\mathbf{h}}$ intersects $R$ non-trivially. Hence $R$
  does not admit a complement in $L$.  \hfill $\Diamond$
\end{eje}

We translate our results for $\Z_p$-Lie lattices into corresponding
statements about $p$-adic analytic pro-$p$ groups.

\begin{pro} \label{solvable_radical_group} Let $G$ be a $p$-adic
  analytic pro-$p$ group. Then $G$ has a unique maximal soluble normal
  subgroup $R$. Moreover, if $G$ is saturable, then $R$ is isolated in
  $G$, i.e.\ $R = \iso_G(R)$, and $R$ is PF-embedded in $G$.
\end{pro}

\begin{proof}
  As $G$ has finite rank, every subgroup of $G$ is finitely
  generated. This implies that $G$ has a unique maximal soluble normal
  subgroup.

  Suppose that $G$ is saturable, and consider the $\Z_p$-Lie lattice
  $L := L(G)$ associated to $G$. Recall that $G$ and $L$ are the same
  as sets. Let $R$ be the soluble radical of the saturable Lie lattice
  $L$. By Proposition~\ref{solvable_radical}, $R$ is isolated and
  PF-embedded in $L$. In particular, $R$ is saturable. According to
  \cite[Theorem 4.5 and Corollary 4.7]{GS}, the subset $R$ forms a
  soluble PF-embedded subgroup of the saturable group $G$. Clearly,
  $R$ is also isolated as a group.

  We claim that $R$ is the maximal soluble normal subgroup of $G$. Let
  $N$ be a soluble normal subgroup of $G$. We have to show that $N
  \subseteq R$. By Proposition~\ref{isolator} we may assume without
  loss of generality that $N$ is saturable. By \cite[\S 4]{GS}, the
  subset $N$ forms a saturable soluble Lie ideal of the $\Z_p$-Lie
  lattice $L$. This shows that $N \subseteq R$.
\end{proof}

\begin{pro} \label{structGroup} Let $G$ be a saturable pro-$p$ group
  with soluble radical $R$. Then $G$ is virtually isomorphic to a
  semidirect product of $H$ and $R$, where $H$ is the direct product
  of saturable pro-$p$ groups whose corresponding $\Q_p$-Lie algebras
  are simple.
\end{pro}

\begin{proof}
  Consider the $\Q_p$-Lie algebra $\mathcal{L} := \mathcal{L}(G)$. Let
  $\mathcal{R}$ denote the soluble radical of $\mathcal{L}$. Then
  there exists a semisimple subalgebra $\mathcal{H}$ of $\mathcal{L}$
  such that $\mathcal{L}$ is the semidirect sum of $\mathcal{H}$ and
  $\mathcal{R}$. As $\mathcal{H}$ is semisimple it decomposes as a
  direct sum $\mathcal{H} = \mathcal{H}_1 \oplus \ldots \oplus
  \mathcal{H}_m$ of simple subalgebras.

  Now consider the $\Z_p$-Lie lattice $L := L(G)$ associated to
  $G$. We know that $R := L \cap \mathcal{R}$, the soluble radical of
  $L$, is isolated in $L$ and hence saturable. For $i \in
  \{1,\ldots,m\}$ put $H_i := L \cap \mathcal{H}_i$. Then each $H_i$
  is isolated in $L$ and hence saturable. From the characterisation of
  saturable $\Z_p$-Lie lattices in terms of potent filtrations (cf.\
  the remark following Theorem~\ref{equivalencia}) it follows that the
  direct sum $H := H_1 \oplus \ldots \oplus H_m$ is saturable. Observe
  that $H \cap R = \{0\}$ and that the semidirect sum $H + R$ is open
  in $L$. The proposition now follows by passing to the saturable
  subgroups of $G$ corresponding to $H$ and $R$.
\end{proof}

\begin{cor}
  \label{saturablejustinfinite} Let $G$ be a saturable insoluble
  just-infinite pro-$p$ group. Then the associated $\Q_p$-Lie algebra
  $\mathcal{L}(G)$ is simple, and $G$ is hereditarily just-infinite.
\end{cor}

\begin{proof}
  For a contradiction suppose that $\mathcal{L} := \mathcal{L}(G)$ is
  not simple. Since $G$ has no soluble normal subgroups, the soluble
  radical of $\mathcal{L}$ is trivial. The semisimple Lie algebra
  $\mathcal{L}$ decomposes as a direct sum of simple Lie algebras
  $\mathcal{H}_1 \oplus \ldots \oplus \mathcal{H}_m$ where $m \geq
  2$. Then $H_1 := L \cap \mathcal{H}_1$ is an isolated and hence
  saturable Lie ideal of $L = L(G)$. It corresponds to a non-trivial
  normal subgroup of infinite index in $G$. This is the required
  contradiction.
\end{proof}


\section{Just-infinite groups}\label{sec:just-infinite}

As indicated in the introduction, it is natural to ask whether
Proposition~\ref{p-1} can be strengthened in a more restricted
setting. Of particular interest in this context is the class of
$p$-adic analytic just-infinite pro-$p$ groups. These groups, also
termed $p$-adically simple groups, are studied to great depth in
\cite{KlLePl97}. In the present section we explore what implications our
approach has within this more restricted class of $p$-adic analytic
pro-$p$ groups. There are two basic questions: when are these groups
torsion-free and when are they saturable? In view of
Theorem~\ref{teo_A} we are particularly interested in answering these
questions in dimension equal to or slightly larger than $p$.

There are two types of $p$-adic analytic just-infinite pro-$p$ groups,
soluble and insoluble. The soluble ones are irreducible $p$-adic space
groups and have dimension $(p-1)p^r$, $r \in \N_0$;
see~\cite[Section~10]{LeMc02}. Hence the smallest dimension above
$p-1$ is $(p-1)p$; in particular, dimension $p$ occurs in this context
only for $p = 2$. But more significantly, soluble just-infinite
pro-$p$ groups simply fail to be torsion-free.

\begin{pro}
  Every soluble just-infinite pro-$p$ group other than $\Z_p$ has
  torsion.
\end{pro}

\begin{proof}
  Let $G$ be a soluble just-infinite pro-$p$ group. Then $G$ is
  virtually abelian. Let $A$ be a maximal abelian normal subgroup of
  $G$. Then $A$ is torsion-free and self-centralising in $G$. If $G =
  A$, then $G \cong \Z_p$ and there is nothing further to show. Now
  assume that $A$ is properly contained in $G$ and write $H :=
  G/A$. Then $V := \Q_p \otimes A$ is a faithful $H$-module. Choose $g
  \in G$ such that $h := gA$ is central and of order $p$ in $H$. By
  Clifford's theorem $V$ decomposes into a direct sum of irreducible
  $\langle h \rangle$-submodules of equal dimensions. Since $\langle h
  \rangle$ acts faithfully on $V$, the trivial representation does not
  occur in this decomposition. But clearly $h$ fixes $g^p \in A$, and
  therefore $g^p = 1$.
\end{proof}

The insoluble $p$-adic analytic just-infinite pro-$p$ groups can be
realised as open subgroups of the groups of $\Q_p$-rational points of
certain semisimple algebraic groups;
cf.~\cite[Section~III]{KlLePl97}. The components of the corresponding
semisimple Lie algebras are pairwise isomorphic and their number is a
power of $p$. Every simple Lie algebra $\mathcal{L}$ over $\Q_p$ is
absolutely simple over its centroid
$\textup{End}_{\mathcal{L}}(\mathcal{L})$. The dimension of an
insoluble $p$-adic analytic just-infinite pro-$p$ group is therefore
of the form $p^k m d$, where $p^k$ is the number of simple components,
$m$ is the dimension of the centroid and $d$ is the dimension of an
absolutely simple Lie algebra. From the classification of absolutely
simple Lie algebras we conclude that dimension $p$ occurs only for
$p=3$.

We note that conversely, every semisimple $\Q_p$-Lie algebra with
simple components which are pairwise isomorphic and whose number is a
power of $p$ occurs as the Lie algebra associated to some insoluble
$p$-adic analytic just-infinite pro-$p$ group;
cf.~\cite[Proposition~(III.9)]{KlLePl97}. We deduce a slightly
stronger result of relevance in our context.

\begin{pro}
\label{justinfinitenosimple}
Let $\mathcal{L}$ be a simple $\Q_p$-Lie algebra and let $k \in
\N_0$. Then there exists a torsion-free insoluble $p$-adic analytic
just-infinite pro-$p$ group $G$ such that the $\Q_p$-Lie algebra
associated to $G$ is isomorphic to the direct sum of $p^k$ copies of
$\mathcal{L}$.
\end{pro}

\begin{proof}
  Let $H$ be a torsion-free open pro-$p$ subgroup of
  $\textup{Aut}(\mathcal{L})$ and consider the standard wreath product
  $W := H \wr C$ where $C = \langle x \rangle$ is a cyclic group of
  order $p^k$. According to \cite[Proposition~(III.9)]{KlLePl97} it
  suffices to find a torsion-free open subgroup of $W$ which projects
  onto $C$ under the natural projection.

  Choose a non-trivial element $h \in H$ and an open normal subgroup
  $N \trianglelefteq H$ such that $h^{p^k} \not \in N$. The base group
  of $W$ is the $p^k$-fold product $H \times \ldots \times H$ and
  accordingly we write its elements as $p^k$-tuples in $H$. Let $G$ be
  the group generated by $x(h,\ldots,h)$ and the open normal subgroup
  $N \times \ldots \times N \trianglelefteq W$. Clearly, $G$ is an
  open subgroup of $W$ projecting onto $C$. It remains to prove that
  $G$ is torsion-free. As $H$ is torsion-free this follows from the
  observation that
  $$
  (x(h,\ldots,h))^{p^k} = (h^{p^k},\ldots,h^{p^k}) \in (H \times
  \ldots \times H) \setminus (N \times \ldots \times N).
  $$
\end{proof}

In \cite{Kl05a} it was shown by example that torsion-free insoluble
$p$-adic analytic maximal just-infinite pro-$p$ groups need not be
saturable. Here we provide further examples:
Corollary~\ref{saturablejustinfinite} shows that the groups
constructed in Proposition~\ref{justinfinitenosimple} for $k \geq 1$
are not saturable. As indicated in the introduction, torsion-free
$p$-adic analytic just-infinite pro-$p$ groups whose $\Q_p$-Lie
algebras are isomorphic to a sum of $p$ copies of
$\mathfrak{sl}_2(\Q_p)$ yield examples of $3p$-dimensional groups
which are not saturable.

Moreover, combining Corollary~\ref{saturablejustinfinite} and
Proposition~\ref{justinfinitenosimple} we obtain
Proposition~\ref{pro_F}.  Another interesting problem concerning
saturable just-infinite pro-$p$ groups is recorded in
Question~\ref{que:all_open_saturable}.


\section{Pro-$p$ groups of dimension at most
  $3$} \label{sec:classification}

Taking advantage of Lazard's correspondence, finite $p$-groups of
reasonably small order can be classified or at least enumerated by
solving the respective problem for finite nilpotent Lie rings of
$p$-power order; e.g.\ see~\cite{OBVa05}. The point is that, because
of their linear structure, Lie rings are generally more tractable than
groups. Being based on the exponential and logarithm maps, Lazard's
correspondence is restricted to groups and Lie rings of nilpotency
class less than $p$.

As an illustration of the method, let us consider non-abelian groups
and nilpotent Lie rings of order $p^3$. One easily verifies that there
are precisely two isomorphism types of nilpotent Lie rings of order
$p^3$, represented by
\begin{align*}
  L_1 & = \langle x,y \rangle_+ & & \text{where $px = p^2 y = 0$ and
    $[y,x] = py$,} \\
  L_2 & = \langle x,y,z \rangle_+ & & \text{where $px = py = pz = 0$
    and $[x,y] = z$ is central.}
\end{align*}
Applying Lazard's correspondence, for $p \geq 3$ we deduce that there
are precisely two non-abelian $p$-groups of order $p^3$, namely
\begin{align*}
  G_1 & = \langle x,y \mid x^p = y^{p^2} = 1 \text{ and } [x,y] = y^p
  \rangle, \\
  G_2 & = \langle x,y \mid x^p = y^p = [x,y]^p = 1 \text{ and } [x,y]
  \text{ central}\rangle.
\end{align*}

In the present section we show how a similar approach, based on
Theorem~\ref{teo_B}, leads to a method for classifying torsion-free
$p$-adic analytic pro-$p$ groups of dimension less than $p$. We work
out the details only in the soluble case up to dimension $3$. It goes
without saying that, while illustrating the point, our results in this
direction are much less sophisticated than the work in \cite{OBVa05},
say. It can be expected that for groups of dimension $4$ and higher
the classification problem becomes rather difficult.

Clearly, there is only one torsion-free $p$-adic analytic pro-$p$
group of dimension $1$, namely the infinite procyclic group $\langle x
\rangle \cong \Z_p$. More generally there is precisely one abelian
torsion-free $p$-adic analytic pro-$p$ group in every given dimension.
Our task is to construct the non-abelian groups in dimensions $2$ and
$3$ from the corresponding Lie lattices.

\subsection{Groups of dimension $2$}
Consider a non-abelian $\Z_p$-Lie lattice $L$ of dimension
$2$. Tensoring with $\Q_p$, we obtain a non-abelian Lie algebra $\Q_p
\otimes L $ of dimension $2$. This Lie algebra can be spanned by
elements $x,y$ with the Lie product given by $[y,x] = y$.

From this we can describe all non-abelian $\Z_p$-Lie lattices of
dimension $2$. Every Lie lattice $L$ of this type can be spanned by
elements $x,y$ whose Lie product is given by $[y,x] = p^s y$ with $s
\in \N_0$. Moreover, the invariant $s = s(L)$ distinguishes the
isomorphism class of $L$, as it determines $|C_L([L,L]):[L,L]| = p^s$,
the index of the commutator Lie sublattice in its centraliser. Note
that $L$ is residually-nilpotent if and only if
$s > 0$. Theorem~\ref{teo_B} allows us to translate our findings to
groups.

\begin{pro}
  For $p \geq 3$ there exists precisely one infinite family of
  $2$-dimensional non-abelian torsion-free $p$-adic analytic pro-$p$
  groups, parameterised by $s \in \N$, namely
  $$
  G(s) = \langle x,y \mid [x,y] = y^{p^s} \rangle.
  $$
\end{pro}

We remark that these groups are not only saturable, but even uniformly
powerful; cf.~\cite[Section~7.3]{Kl03}. Furthermore we have $G(s+1)
\cong G(s)^p$ for all $s \in \N$, and the $2$-dimensional abelian
torsion-free $p$-adic analytic pro-$p$ group can be regarded as a
limit of the family $G(s)$ as $s$ tends to infinity.

For completeness we also consider the case $p=2$, where the Lie
correspondence breaks down. It is convenient to define $2^\infty := 0$
in $\Z_2$.

\begin{pro}
  There exist precisely two infinite families of $2$-dimen\-sional
  torsion-free $2$-adic analytic pro-$2$ groups, each parameterised by
  $s \in \N \cup \{ \infty \}$ with $s \geq 2$, namely
  $$
  G_+(s) = \langle x,y \mid [x,y] = y^{2^s} \rangle \quad \text{and}
  \quad G_-(s) = \langle x,y \mid [x,y] = y^{-2-2^s} \rangle.
  $$
\end{pro}

\begin{proof}
  The cyclic subgroups of $\textup{Aut}(\Z_2) \cong \Z_2^* = \langle
  -1 \rangle \times \langle 1 + 4 \rangle$ are precisely the groups
  $\langle 1 + 2^s \rangle$ and $\langle -1-2^s \rangle$, where $s \in
  \N \cup \{\infty\}$ with $s \geq 2$. Hence it suffices to show that
  every torsion-free $2$-adic analytic pro-$2$ group is isomorphic to
  a semidirect product of $\Z_2$ by $\Z_2$.

  Let $G$ be a torsion-free $2$-adic analytic pro-$2$ group. Then $G$
  contains a saturable open subgroup $H$. By consideration of the
  associated $\Z_2$-Lie lattice one sees that $H$ is of the form $\Z_2
  \ltimes \Z_2$, hence isomorphic to one of the groups $G_+(s)$,
  $G_-(s)$ for suitable $s \in \N \cup \{\infty\}$ with $s \geq
  2$. Arguing by induction on $\lvert G:H \rvert$, we may assume that
  $\lvert G:H \rvert = 2$.

  \smallskip

  \noindent \textit{Case 1:} $H$ is abelian. Write $H = \Z_2 x + \Z_2
  y \cong \Z_2 \times \Z_2$ and $G = \langle g \rangle H$. If $g$ acts
  trivially on $H$, then $G$ is abelian and there is nothing further
  to prove. Hence suppose that $g$ acts on $H$ as an involution. The
  finite $2$-subgroups of $\textup{GL}_2(\Z_2)$ are known;
  cf.~\cite[Section~10]{LeMc02}. Carrying out if necessary a change of
  basis, we may assume that the action of $g$ on $H$ with respect to
  the basis $(x,y)$ is given by one of the following pairs of
  equations
  \begin{enumerate}
  \item[(a)] $x^g = x^{-1}$ and $y^g = y$,
  \item[(b)] $x^g = x^{-1}$ and $y^g = y^{-1}$,
  \item[(c)] $x^g = y$ and $y^g = x$.
  \end{enumerate}
  Suppose that $g$ acts on $H$ according to the equations (a). Then
  $g^2 \in Z(G) = \langle y \rangle$, and since $G$ is torsion-free,
  $g^2$ must be a generator of the procyclic group $\langle y
  \rangle$. This shows that $G = \langle g \rangle \ltimes \langle x
  \rangle$ is isomorphic to a semidirect product of $\Z_2$ by $\Z_2$.

  The second and third type of action do not actually arise. If the
  action of $g$ on $H$ was given by the equations (b), then $g^2 \in
  Z(G) = 1$ and $G$ would not be torsion-free.

  Now assume for a contradiction that $g$ acts on $H$ according to the
  equations (c). Then $g^2 \in Z(G) = \langle xy \rangle$, and since
  $G$ is torsion-free, $g^2$ must be a generator of the procyclic
  group $\langle xy \rangle$. Without loss of generality, we may
  assume that $g^2 = xy$. Then $gy^{-1} \not = 1$, but $(gy^{-1})^2 =
  g^2 x^{-1} y^{-1} = 1$, in contradiction to $G$ being torsion-free.

  \smallskip

  \noindent \textit{Case 2:} $H$ is non-abelian. In this case $H =
  \langle x \rangle \ltimes \langle y \rangle$ is a split extension,
  where
  \begin{equation} \label{eq:conj_of_y_by_x} x^{-1} y x = y^{\pm (1 +
      2^s)} \not = y \quad \text{for suitable $s \in \N \cup
      \{\infty\}$ with $s \geq 2$.}
  \end{equation}
  We observe that the $1$-dimensional subgroup $N := \langle y \rangle
  = \iso_H([H,H])$ is characteristic in $H$ and hence normal in
  $G$. The quotient $G/N$ is, of course, $1$-dimensional. In
  particular, if $G/N$ is torsion-free, then the exact sequence $1
  \rightarrow N \rightarrow G \rightarrow G/N \cong \Z_2 \rightarrow
  1$ splits, and $G$ is isomorphic to a semidirect product of $\Z_2$
  and $\Z_2$. Similarly, if $G/N$ is abelian, then writing $M :=
  \iso_G(N) \cong \Z_2$, the exact sequence $1 \rightarrow M
  \rightarrow G \rightarrow G/M \cong \Z_2 \rightarrow 1$ splits, and
  $G$ is isomorphic to a semidirect product of $\Z_2$ and $\Z_2$.

  Hence we assume, for a contradiction, that $G/N$ is non-abelian and
  has torsion. Since $H/N$ is torsion-free, we find $g \in G$ with $G
  = \langle g \rangle H$ and $g^2 \in N$. As $G$ is torsion-free, the
  $1$-dimensional group $\langle g \rangle N$ is procyclic and $g^2$
  must be a generator of $N = \langle y \rangle$. Without loss of
  generality we may assume that $g^2 = y$. Being non-abelian, $G/N =
  \langle gN \rangle \ltimes \langle xN \rangle$ is isomorphic to the
  infinite pro-$2$ dihedral group. Hence we find $k \in \Z$ such that
  $g^{-1} x g = x^{-1} y^k = x^{-1} g^{2k}$. Writing $w := g^{-1} x
  \not = 1$, this gives $w^2 = g^{2(k-1)} = y^{k-1}$.

  By~\eqref{eq:conj_of_y_by_x} this implies
  $$
  w^2 = w^{-1} g^{2(k-1)} w = x^{-1} y^{k-1} x = y^{\pm (1+2^s)(k-1)}
  = (w^2)^{\pm (1+2^s)}
  $$
  where $\pm (1+2^s) \not = 1$. Hence $w$ has finite order, in
  contradiction to $G$ being torsion-free.
\end{proof}


\subsection{Soluble groups of dimension $3$}
Again we aim for a classification of groups by first describing the
relevant Lie lattices. Since the correspondence between saturable
$\Z_p$-Lie lattices and saturable pro-$p$ groups breaks down for $p
\leq 3$ in dimension $3$, we avoid from the beginning the extra work
required to deal with the prime $2$: throughout the present
subsection, let $p$ denote an odd prime.

Consider a soluble $\Z_p$-Lie lattice $L$ of dimension $3$. Tensoring
with $\Q_p$, we obtain a soluble Lie algebra $\mathcal{L} := \Q_p
\otimes L$ of dimension $3$. According to \cite[Chapter 1~\S4]{Ja},
the Lie algebra $\mathcal{L}$ has an abelian ideal of dimension
$2$. Consequently $\mathcal{L}$ decomposes as the semidirect sum of a
$1$-dimensional subalgebra and such a $2$-dimensional abelian
ideal. If $[\mathcal{L},\mathcal{L}]$ is central, then $\mathcal{L}$
is nilpotent and consequently either abelian or isomorphic to the
so-called Heisenberg Lie algebra. If $[\mathcal{L},\mathcal{L}]$ is
not central, then $\mathcal{L}$ fails to be nilpotent and the
centraliser $C_{\mathcal{L}}([\mathcal{L},\mathcal{L}])$ is the unique
$2$-dimensional abelian ideal of $\mathcal{L}$.

Our discussion shows that, in any case, the original Lie lattice $L$
has an abelian ideal $I$ of dimension $2$ such that $L/I$ is a Lie
lattice of dimension $1$. Write $I = \Z_p y_1 + \Z_p y_2$, and choose
an additive complement $K = \Z_p x$ to $I$ in $L$. Then $L = \Z_p x +
\Z_p y_1 + \Z_p y_2$ decomposes as a semidirect sum $K \ltimes I$; the
Lie bracket is determined by the equations
$$
[y_1,y_2] = 0, \quad [y_1,x] = a_{11} y_1 + a_{12} y_2, \quad [y_2,x]
= a_{21} y_1 + a_{22} y_2
$$
for suitable coefficients $a_{ij} \in \Z_p$. Define
$$
A :=
\begin{pmatrix} a_{11} & a_{12} \\ a_{21} & a_{22}
\end{pmatrix}\in \mathfrak{gl}_2(\Z_p).
$$
The Lie lattice $L$ is abelian if and only if $A$ is the null
matrix. Similarly, it is nilpotent if and only if $A$ is
nilpotent. These cases are fairly easy to deal with. On the other
hand, if $L$ is not nilpotent, then $I = C_L([L,L])$ is in fact
uniquely determined. The matrix $A$, however, generally depends on the
particular basis $(x,y_1,y_2)$. Changing this basis corresponds to
conjugating and scaling the matrix $A$. This introduces an equivalence
relation on the set $\mathfrak{gl}_2(\Z_p)$: the
\emph{multiplicative-similarity class} of $A$,
$$
\mathbf{A}_L := \{ u B^{-1} A B \mid u \in \Z_p^* \text{ and } B \in
\textup{GL}_2(\Z_p) \},
$$
does not depend on the chosen basis $(x,y_1,y_2)$. Moreover, by our
observations the invariant $\mathbf{A}_L$ characterises and
parameterises the isomorphism class of the Lie lattice $L$.

\begin{pro}[Representatives for multiplicative-similarity
  classes] \label{list} Every nilpotent matrix $A \in
  \mathfrak{gl}_2(\Z_p)$ is multiplicatively-similar to precisely one
  of the matrices
  $$
  \begin{pmatrix} 0 & 0 \\ 0 & 0
  \end{pmatrix}
  \qquad \text{or} \qquad p^s \begin{pmatrix} 0 & 0 \\ 1 &
    0 \end{pmatrix} \text{, $s \in \N_0$.}
  $$

  Let $\rho \in \Z_p^*$, not a square modulo $p$. Then every
  non-nilpotent matrix $A \in \mathfrak{gl}_2(\Z_p)$ is
  multiplicatively-similar to precisely one matrix of the form $p^s
  A_0$, where $s \in \N_0$ and $A_0$ is one of the following core
  matrices:
  \begin{enumerate}
  \item $\begin{pmatrix} 1 & 0 \\ 0 & 1 \end{pmatrix}$;
  \item $\begin{pmatrix} 1 & 0 \\ 0 & 1 \end{pmatrix} +
    p^r \begin{pmatrix} 0 & d \\ 1 & 0 \end{pmatrix}$, where $r \in
    \N$ and $d \in \Z_p$;
  \item $\begin{pmatrix} 0 & d \\ 1 & p^r \end{pmatrix}$, where $r
    \in \N_0$ and $d \in \Z_p$;
  \item $\begin{pmatrix} 0 & p^r \\ 1 & 0 \end{pmatrix}$ or
    $\begin{pmatrix} 0 & \rho p^r \\ 1 & 0 \end{pmatrix}$, where $r
    \in \N_0$.
  \end{enumerate}
\end{pro}

\begin{proof}
  Recall that we assume $p > 2$ throughout. From the analysis in
  \cite{ApOn83}, we deduce that the $\textup{GL}_2(\Z_p)$-conjugacy
  class of a matrix $A \in \mathfrak{gl}_2(\Z_p)$ which is not scalar
  modulo $p$ is uniquely determined by its characteristic polynomial,
  i.e.\ by its trace and determinant.

  Clearly, the null matrix forms a multiplicative-similarity class by
  itself. Thus it suffices to consider a non-zero matrix $A \in
  \mathfrak{gl}_2(\Z_p)$ and its multiplicative-similarity class
  $\mathbf{A}$. Write $A = p^s A_0$ where $s \in \N_0$ is chosen so
  that $A_0 \not \equiv 0$ modulo $p$. Note that $s$ is in fact an
  invariant of the class $\mathbf{A}$, independent of the particular
  representative $A$.

  It is enough to concentrate on $A_0$ and its
  multiplicative-similarity class $\mathbf{A}_0$. We distinguish four
  cases and show that $A_0$ is multiplicatively-similar to one of the
  listed core matrices. At the same time it will become clear that the
  listed matrices are pairwise not multiplicatively-similar.

  \smallskip

  \noindent \textit{Case 1:} $A_0$ is scalar. Multiplying by a
  suitable unit $u \in \Z_p^*$ (uniquely determined), we see that
  $A_0$ is multiplicatively-similar to $\left(\begin{smallmatrix} 1 &
      0 \\ 0 & 1 \end{smallmatrix} \right)$.

  \smallskip

  \noindent \textit{Case 2:} $A_0$ is scalar modulo $p$, but not
  scalar. Let $r \in \N$ such that $A_0$ is scalar modulo $p^r$, but
  not scalar modulo $p^{r+1}$. Note that $r$ is in fact an invariant
  of the class $\mathbf{A}_0$, independent of the particular
  representative $A_0$. Multiplying by a suitable unit $u_1 \in
  \Z_p^*$ (uniquely determined modulo $p^r$) we see that $A_0$ is
  multiplicatively-similar to $\left( \begin{smallmatrix} 1 & 0 \\ 0 &
      1 \end{smallmatrix} \right) + p^r A_1$. Multiplying by a
  suitable one-unit $u_2 \in 1 + p^r \Z_p$ (uniquely determined), we
  see that $A_0$ is multiplicatively-similar to
  $\left( \begin{smallmatrix} 1 & 0 \\ 0 & 1 \end{smallmatrix} \right)
  + p^r A_2$ where $A_2$ has trace $0$. Define $d := -\det(A_2) \in
  \Z_p$. Conjugating by a suitable element of $\textup{GL}_2(\Z_p)$,
  we see that $A_0$ is multiplicatively-similar to
  $\left( \begin{smallmatrix} 1 & 0 \\ 0 & 1 \end{smallmatrix} \right)
  + p^r \left( \begin{smallmatrix} 0 & d \\ 1 & 0 \end{smallmatrix}
  \right)$.

  \smallskip

  \noindent \textit{Case 3:} $A_0$ is not scalar modulo $p$ and has
  non-zero trace. Scaling by a unit $u \in \Z_p^*$ (uniquely
  determined), we find that $A_1 := uA_0$ has trace $p^r$ for suitable
  $r \in \N_0$. Note that $r$ is in fact an invariant of the class
  $\mathbf{A}_0$, independent of the particular representative
  $A_0$. Put $d := - \det(A_2)$. Conjugating by a suitable element of
  $\textup{GL}_2(\Z_p)$, we see that $A_0$ is multiplicatively-similar to
  $\left( \begin{smallmatrix} 0 & d \\ 1 & p^r \end{smallmatrix}
  \right)$.

  \smallskip

  \noindent \textit{Case 4:} $A_0$ is not scalar modulo $p$ and has
  zero trace. Scaling by a unit $u \in \Z_p^*$, we find that $A_1 :=
  uA_0$ retains trace $0$ and satisfies $\det(A_1) = 0$ or $ \det(A_1)
  \in \{-p^r, -\rho p^r \}$ for suitable $r \in \N_0$. Note that in
  the latter case $r$ is an invariant of the class $\mathbf{A}_0$,
  independent of the particular representative $A_0$.  Conjugating by
  a suitable element of $\textup{GL}_2(\Z_p)$, we see that $A_0$ is
  multiplicatively-similar to one of the matrices
  $\left( \begin{smallmatrix} 0 & 0 \\ 1 & 0 \end{smallmatrix}
  \right)$, $\left( \begin{smallmatrix} 0 & p^r \\ 1 &
      0 \end{smallmatrix} \right)$, $\left( \begin{smallmatrix} 0 &
      \rho p^r \\ 1 & 0 \end{smallmatrix} \right)$.
\end{proof}

Let $A$ be one of the matrices listed in Proposition~\ref{list}, and
let $L$ be the $\Z_p$-Lie lattice associated to the
multiplicative-similarity class $\mathbf{A}$ of $A$. Then $L$ is
residually-nilpotent if and only if $A^2 \equiv_p 0$. Under the
hypothesis $p > 3$ Theorem~\ref{teo_B} allows us to translate our
findings to groups. For this it is convenient to embed $L$ as a Lie
sublattice into $\mathfrak{gl}_3(\Z_p)$ by mapping
$$
x \mapsto
\begin{pmatrix}
  A & \begin{matrix} 0 \\ 0 \end{matrix} \\
  \begin{matrix} 0 & 0 \end{matrix} & 0
\end{pmatrix}
\quad y_1 \mapsto \begin{pmatrix} 0 & 0 & 0 \\ 0 & 0 & 0 \\ 1 & 0 &
  0 \end{pmatrix}, \quad y_2 \mapsto \begin{pmatrix} 0 & 0 & 0 \\ 0 &
  0 & 0 \\ 0 & 1 & 0 \end{pmatrix}.
$$
If $L$ is residually-nilpotent, the corresponding group $G$ can be
obtained by applying the ordinary exponential map: $G = H \ltimes N$
is the semidirect product of $H \cong \Z_p$ and a free abelian pro-$p$
group $N \cong \Z_p \times \Z_p$. The action of $H$ on $N$ is given by
$\exp(A) = \sum_{n=0}^\infty A^n/n!$.

Alternatively, we can use Proposition~\ref{list} directly to write
down finite presentations for representatives $G = H \ltimes N$ of
the isomorphism classes of $3$-dimensional soluble torsion-free
$p$-adic analytic pro-$p$ groups. Namely, the action of $H$ on $N$ can
be given by $1+A$ rather than $\exp(A)$. Here $A$ runs again through
those matrices listed in Proposition~\ref{list} which are nilpotent
modulo $p$.

\begin{teo}\label{thm:complete_list}
  For $p > 3$ the following constitutes a complete and irredundant
  list of $3$-dimensional torsion-free $p$-adic analytic pro-$p$
  groups (up to isomorphism):
  \begin{enumerate}
  \item the abelian group, isomorphic to $\Z_p^3$,
    $$
    G_0(\infty) = \langle x_1, x_2, x_3 \mid [x_1,x_2] = [x_1,x_3] =
    [x_2,x_3] = 1 \rangle.
    $$
  \item the two-step nilpotent groups, parameterised by $s \in \N_0$,
    $$
    G_0(s) = \langle x, y, z \mid [x,y] = z^{p^s}, [x,z]=[y,z]=1
    \rangle;
    $$
    these are open subgroups of the Heisenberg group, viz.\ $G_0(0)$.
  \item a family of non-nilpotent groups, parameterised by $s \in \N$,
    $$
    G_1(s) = \langle x, y_1, y_2 \mid [y_1,y_2] = 1, [y_1,x] =
    y_1^{p^s}, [y_2,x] = y_2^{p^s} \rangle.
    $$
  \item a family of non-nilpotent groups, parameterised by $s,r \in \N$
    and $d \in \Z_p$,
    \begin{align*}
      G_2(s,r,d) = \langle x, y_1, y_2 & \mid  [y_1,y_2] = 1, \\
      & [y_1,x] = y_1^{p^s} y_2^{p^{s+r}d}, [y_2,x] = y_1^{p^{s+r}}
      y_2^{p^s} \rangle.
    \end{align*}
  \item a family of non-nilpotent groups, parameterised by $s,r \in
    \N_0$ and $d \in \Z_p$ such that \textup{(a)} $s \geq 1$ or
    \textup{(b)} $r \geq 1$ and $d \in p\Z_p$,
    \begin{align*}
      G_3(s,r,d) = \langle x, y_1, y_2 & \mid [y_1,y_2] = 1, \\
      & [y_1,x] = y_2^{p^s d}, [y_2,x] = y_1^{p^s} y_2^{p^{s+r}}
      \rangle.
    \end{align*}
  \item two families of non-nilpotent groups, parameterised by $s,r \in
    \N_0$ with $s+r \geq 1$,
    \begin{align*}
      G_4(s,r) & = \langle x, y_1, y_2 \mid [y_1,y_2] = 1, [y_1,x] =
      y_2^{p^{s+r}}, [y_2,x] = y_1^{p^s} \rangle,\\
      G_5(s,r) & = \langle x, y_1, y_2 \mid [y_1,y_2] = 1, [y_1,x] =
      y_2^{p^{s+r} \rho}, [y_2,x] = y_1^{p^s} \rangle.\\
    \end{align*}
  \end{enumerate}
\end{teo}

We remark that among these groups there are several which are not
uniformly powerful. These complement the insoluble examples of
saturable but not uniformly powerful pro-$p$ groups given in
\cite{Kl05a}.


\subsection{Insoluble groups of dimension $3$}

Tensoring an insoluble $\Z_p$-Lie lattice $L$ with $\Q_p$, we obtain
an insoluble Lie algebra $\mathcal{L} := \Q_p \otimes L$. In dimension
$3$ the Lie algebra $\mathcal{L}$ is necessarily simple of type
$\textup{A}_1$ and only two forms occur. They are represented by
$\mathfrak{sl}_2(\Q_p)$ and $\mathfrak{sl}_1(\mathbb{D}_p)$, where
$\mathbb{D}_p$ denotes a central division algebra of index $2$ over
$\Q_p$.

Suppose that $p > 3$, so that Theorems~\ref{teo_A} and \ref{teo_B}
become applicable. Any saturable $3$-dimensional insoluble pro-$p$
group $G$ acts faithfully on itself, regarded as a Lie lattice
$L$. This adjoint action provides an embedding of $G$ as an open
subgroup into $\textup{Aut}(\mathcal{L})$ where $\mathcal{L} = \Q_p
\otimes L$ as before.

It is known that the open maximal pro-$p$ subgroups of
$\textup{Aut}(\mathcal{L})$ form a unique conjugacy class;
see~\cite[Lemma~(III.16)]{KlLePl97}. Any open maximal pro-$p$ subgroup
of $\textup{Aut}(\mathfrak{sl}_2(\Q_p))$ is isomorphic to the group
$$
\textup{SL}_2^{\vartriangle}(\Z_p) := \{ g \in \textup{SL}_2(\Z_p)
\mid \text{$g$ upper uni-triangular modulo $p$} \},
$$
a Sylow pro-$p$ subgroup of $\textup{SL}_2(\Z_p)$;
cf.~\cite[Lemma~(XI.4)]{KlLePl97}. Denote the maximal $\Z_p$-order of
$\mathbb{D}_p$ by $\Delta_p$. Then any open maximal pro-$p$ subgroup
of $\textup{Aut}(\mathfrak{sl}_1(\mathbb{D}_p))$ is isomorphic to
$\textup{SL}_1^1(\Delta_p)$, a Sylow pro-$p$ subgroup of
$\textup{SL}_1(\mathbb{D}_p)$; see \cite[Proof of
Proposition~10.4]{Kl03}. For $p > 3$, both of these groups are
torsion-free and hence saturable; cf.~\cite[Theorem~1.3]{Kl05a}.

In view of Theorem~\ref{teo_A}, the $3$-dimensional insoluble
torsion-free $p$-adic analytic pro-$p$ groups are precisely the open
subgroups of $\textup{SL}_2^{\vartriangle}(\Z_p)$ and
$\textup{SL}_1^1(\Delta_p)$. By Theorem~\ref{teo_B} they correspond to
open Lie sublattices of the corresponding $\Z_p$-Lie lattices
\begin{align*}
  \mathfrak{sl}_2^{\vartriangle}(\Z_p) & \cong \langle x,y,h \mid
  [x,y] = h, [x,h] = -2px, [y,h] = 2py \rangle, \\
  \mathfrak{sl}_1^1(\Delta_p) & \cong \langle x,y,z \mid [x,y] = pz,
  [x,z] = p \rho y, [y,z] = -x \rangle.
\end{align*}
Here $\rho \in \Z_p^*$ is not a square modulo $p$. We refrain from
giving an explicit classification up to isomorphism, but the following
fact is perhaps worth pointing out
(cf.~\cite[Lemma~(III.11)]{KlLePl97}): an insoluble just-infinite
$p$-adic analytic pro-$p$ group is never isomorphic to one of its
proper subgroups. An interesting treatment of pro-$p$ subgroups of
$\textup{SL}_2(\Z_p)$, valid for $p>2$, can be found in \cite{Pi93}.

\end{document}